%% file: dennisyoung3.tex
\magnification=\magstep1
\input amstex
\input fixitem.tex
\alwaysindent=15pt	
\parindent=0pt
\UseAMSsymbols

\baselineskip = 13pt

\voffset=.5in
\pageheight{7.4in}

\pagewidth{6.45truein}

\footline{\hss\tenrm  -- \folio --\hss}


\def\bsk{\bigskip}
\def\msk{\medskip}
\def\ssk{\smallskip}
\def\cl{\centerline}

\fontdimen7\textfont0=5pt
\def\footnote{\plainfootnote}

\def\bull{{\raise .3ex \hbox{$\ssize \bullet$}}}

\vphantom{}


\def\zp{{\Bbb Z}_p}
\def\slp{\text{SL}_2({\Bbb Z}_p)}
\def\gen#1{\langle #1 \rangle}
\def\ar{\longrightarrow}

\def\char{\vartriangleleft\kern -.7em\vert} 
\def\tr{\operatorname{tr}}
\def\leg#1#2{\fracwithdelims(){#1}{#2}}

\UseAMSsymbols
\NoBlackBoxes
\mathsurround=2pt
\fontdimen7\textfont0=5pt
\fontdimen16\tensy=2.7pt
\fontdimen17\tensy=2.7pt 

\cl{\bf Computation of the Scharlau Invariant, I}
\bsk
\cl{R.~Keith Dennis and Paul K. Young}
\bsk
\bsk

\centerline{\bf Abstract.}
\ssk

The Scharlau invariant determines whether or not a finite group has a
fixed point free representation over a field:\ \ if $0$, yes,
otherwise, no.  Until now it was known to be one of $0$, $1$, $p$,
$p^2$ for $p$ a prime dividing the order of the group.  We eliminate
$p^2$ as a possibility.  Work of Scharlau [Sch] reduces the question
to the above list with $p^2$ being possible for the groups
$\text{SL}_2({\Bbb Z}_p)$ for $p$ a Fermat prime larger than $5$.  A
computation using GAP in the Senior Thesis [Y] of the second author
solves the problem for $p = 17$.  With this motivation, we found a
short proof of the result not requiring a computer.

\footnote{ }{2020 Mathematics Subject Classification. 20C15, 20C05, 20D99.}

\footnote{ }{Key words and phrases. Fixed point free representation,
group ring, Fermat prime.}
 
\msk
\centerline{\bf Introduction.}
\ssk
{\it
\centerline{Insert by first author.}

This short proof was originally written in 2003.  It appears only now
due to the procrastination of the first author alone.  A second part
will add that the problem has been solved, or almost solved, by several
different authors who were not aware of the  others.  In particular,
they did not seem to know the original question of Scharlau.} 
\msk
In the spring of 1975 there was a seminar at Cornell motivated by a
question of George Cooke ``What consequences must always follow from
the action of a finite group on an algebraic structure?''  The main
lectures given by Ken Brown [Br] dealt with the work of Hsia and
Peterson [HP1] and related it to earlier work of Zassenhaus [Z] and
others (see e.g\. [W] Part III).  This note contains the final step in
the determination of the Scharlau invariant which had been defined in
1974 by W.~Scharlau [Sch] to study the question.
\msk
The solution of the problem was made possible by computer computations
of the second author, an undergraduate working on a VIGRE-supported
summer project advised by the first author.  Although he had
substantial prior computing experience, the second author had only just
finished the Honors Abstract Algebra course at Cornell.  This work
makes it clear that with the power available in GAP, it is possible
for beginning students to make contributions to research.
\msk
The problem had been earlier considered in a number of papers (e.g.,
[HP1], [HP2], [K], [B\"o]) and had been raised independently in 1992
by Guralnick and Wiegand [GW].  Their paper in fact included
essentially all previously known results.  The final step involves
computing the Scharlau invariant for the special groups
${\text{SL}_2({\Bbb Z}_r)}$ for $r$ a Fermat prime; the
reduction to that case is essentially a classification theorem of
Suzuki [Su].  The solution involves discovering special invariant
elements in the integral group ring analogous to those constructed by
Scharlau in the case of groups of order $pq$, and even earlier by
Burnside. 
\bsk
\centerline{\bf Short version of proof.}
\msk
Let $\slp$ be the special linear group where $p$ is a prime so that
$p \equiv 17 \mod 60$.  We will later restrict to the case where $p$ is a
Fermat prime which is at least $17$.

\msk
We begin by noting a result going back to I\. Schur [S].
\msk
{\bf Lemma 1.}  $\slp$ has $p+4$ conjugacy classes.
\item{a.}  If the trace is not $\pm 2$, two elements
of $\slp$ are conjugate if and only if their traces are equal,
\item{b.}  If the trace is $2$, there are $3$ conjugacy
classes: \ $I$, $T$ and $Q$.
\item{c.}  If the trace is $-2$, there are $3$ conjugacy classes:
$-I$, $-T$ and $-Q$.
\item{}  Here $T = E_{12}(1)$, and $Q = E_{12}(\alpha)$ are
elementary matrices, and $\alpha \in \zp^*$ is a non-square.

\msk
See Dornhoff [D] Theorem 38.1, page 228 ff for the proof of this and
of Lemma 2.
\msk
{\bf Lemma 2.}
\item{a.} The conjugacy class of a diagonal matrix with trace unequal
to $\pm 2$ has exactly $p(p+1)$ elements.
\item{b.} $\slp$ has an element of order $p+1$.  The conjugacy class
of any element of a subgroup of order $p+1$ with trace unequal to $\pm
2$ has order $p(p-1)$.
\msk
For any finite group $G$, we define a ``centralization'' map $c:\
{\Bbb Z}[G] \ar {\Cal Z}({\Bbb Z}[G])$ as the sum of all distinct
conjugates of $x$:
$$
c(x) = \sum_{g\in G/N} gxg^{-1},
$$
here $N = N(x)$ denotes the centralizer of $x$ in $G$.  One can also
define the ``Maschke'' map $m$ as the sum of all conjugates of $x$.
These two maps are related by $m(x) = |N(x)| c(x)$.
\msk
Let $S$ be a subset of $G$.  We denote by $\overline S$ the element of
${\Bbb Z}[G]$ which is the sum of the elements of $S$.  The {\it
Scharlau ideal} of $G$, is the 2-sided ideal I of ${\Bbb Z}[G]$
generated by the elements $\overline H$ where $H$ ranges over the
non-trivial subgroups of $G$.  This 2-sided ideal is generated over
${\Bbb Z}$ by the elements $\overline{gH}$ for $gH$ ranging over the left
cosets of the non-trivial subgroups of $G$.  In fact it is easy to see
that one need only use the subgroups of prime order.  The {\it
Scharlau invariant} of $G$ is the non-negative integer which generates
$I\cap {\Bbb Z}[G]$.  As a simplification of the notation in the
computations below, the bar over letters will be omitted.  It should
be clear from context if a subset of $G$ or an element of ${\Bbb
Z}[G]$ is under consideration.  For an element $g \in G$ we will write
$\gen{g}$ for the cyclic subgroup generated by $g$.  We write $cl(g)$
for the conjugacy class of the element $g \in G$.
\msk
{\bf Example.}
A partition $P$ of $G$ is a collection of subgroups that have pairwise
trivial intersections and whose union is $G$. Then $|P|-1$ is in 
$I\cap {\Bbb Z}[G]$.
\msk
In case $G = \slp$ Our main result is the following:
\msk
{\bf Theorem.}  The order of $\slp$ is an element of $I$.
\msk
{\bf Corollary.}  The Scharlau invariant of $\slp$ is $p$ for $p$ a Fermat
prime.
\msk
{\bf Proof:} \ It is known that the Scharlau invariant is a non-trivial
divisor of $p^2$ (see e.g\. Guralnick and Wiegand [GW]).  Hence the
greatest common divisor of $p^2$ and $(p-1)p(p+1)$, which is $p$, lies
in $I$ and the result follows.
\msk
We prove the Theorem by giving an explicit formula:
$$
(p-1)p(p+1) = 2(p+1)c(\gen{V}) - 4c(W\gen{V}) + 4(p-1)c(\Delta\gen{T})
\tag{\lower .3 em\hbox{*}}
$$
for certain matrices $V$ and $W$ of order 3 whose product $WV$ is
diagonal.  In fact, the right  hand side of the same equation is 
divisible by $4$ showing that the order of the group over $4$ lies 
in $I$.
\msk
{\bf Remark.}
In case the group $G$ is solvable, the standard examples in the
Scharlau ideal such as for the non-abelian group $G_{pq}$ of order $pq$,
for primes $p$ and $q$, come from a partition with cyclic
subgroups ($c(G_{pq}=p$).  These are invariant under conjugation by $G$.
However, the elements of the Scharlau ideal giving the equation above
are not invariant under conjugation by $G$, contrary to our original
expectations.
\bsk
We do this in several steps.  First a motivational lemma.
\msk
{\bf Lemma 3.}  A conjugacy class of $G$ is a union of left cosets of some
non-trivial subgroup of $G$ if and only if it contains a left coset of
some non-trivial subgroup of $G$.
\msk
The proof is trivial.
\msk
{\bf Lemma 4.}  Let $D$ be any diagonal matrix which is not $\pm I$.  Then
$c(D\gen{T}) = cl(D)$.
\msk
{\bf Proof:} \ The traces of $D$ and $DT^i$ are equal and hence they lie in
the same conjugacy class.  Thus $D\gen{T}$ lies in $cl(D)$ and hence
$cl(D)$ is the union of conjugates of $D\gen{T}$.  Now by Lemma 2
$cl(D)$ has $p(p+1)$ elements, hence $N(D)$ has $p-1$ elements.  As
all diagonal matrices and $T$ lie in $N(D\gen{T})$, it has
order at least $(p-1)p$ and index at most $p+1$.  However the
conjugates of $D\gen{T}$, which contains $p$ elements, must cover all of
$cl(D)$, which has $(p-1)p$ elements, so there must be exactly
$p+1$ such conjugates.  The lemma follows.
\msk
{\bf Lemma 5.}  Let $V$ be the companion matrix of $x^2+x+1$.  Thus $V$ is
an element of order $3$ and trace $-1$.  $c(\gen{V}) = (p-1)p/2 + cl(V)$.
\msk
{\bf Proof:} \ First note that $\gen{V}$ contains $I$, $V$ and
$V^{-1}$ with the latter two conjugate via the matrix
$$
\left(\matrix 0&a\\
a&0\endmatrix\right)
$$
for any $a \in \zp^*$ with $a^2 = -1$.  Such an $a$ exists as $p
\equiv 1 \mod 4$.  Now $|N(V)| = p+1$ by Lemma 2 as the conjugacy
class of $V$ has order $(p-1)p$.  Thus $|N(\gen{V})| = 2(p+1)$.  It
then follows that $c(\gen{V})$ has $(p-1)p(p+1)/2(p+1)$ as the
coefficient of $1$.  There are $3\cdot (p-1)p/2$ terms total in
$c(\gen{V}$ with the remaining $(p-1)p$ of them forming the conjugacy
class of $V$
\msk
{\bf Lemma 6.} \ There exists an element $W \in \slp$ of order $3$
with $WV$ conjugate to $WV^{-1}$.  The latter matrices are conjugate
to a diagonal matrix with trace unequal to $\pm 2$.
\msk
{\bf Proof:} \ Let $W$ be the matrix
$$
\left(\matrix 0&u^{-1}\\
u&-1\endmatrix\right)
$$
which has order $3$ as $\tr W = \tr V = -1$.  We now compute
$WV$
$$
\left(\matrix -u^{-1}&u^{-1}\\
-1&1-u\endmatrix\right)
$$
and $WV^{-1}$
$$
\left(\matrix u^{-1}&0\\
1-u&u\endmatrix\right) \ .
$$
These two matrices will be conjugate precisely when their traces are
equal: 
$$
1 - u - u^{-1} = u + u^{-1}
$$
or
$$
u+u^{-1} = 1/2
$$
or
$$
u^2 -1/2u +1 = 0 \ .
$$
For the known Fermat primes such a $u$ is easily found:
$$
\alignat 4
    p&:\quad&\quad     u&\quad+&\quad u^{-1}&\ \ =\ \ &   1/2&\\
   17&:\quad&\quad     3&\quad+&\quad      6&\ \ =\ \ &     9&\\
  257&:\quad&\quad    28&\quad+&\quad    101&\ \ =\ \ &   129&\\
65537&:\quad&\quad 37162&\quad+&\quad  61144&\ \ =\ \ & 32769&
\endalignat
$$
In the general case, this follows from a computation with Legendre
symbols.  Note that the discriminant of the quadratic equation is
$\sqrt{-15/4\ } = \sqrt{(-1/2)^2-4\ }$.  Thus we need only show
that $\sqrt{-15 } \in \zp^*$.  We have then
$$
\leg{-15}p = \leg{-1}p \leg 3p \leg 5p\ =\ 1 \ :
$$
the first term is $1$ as $p \equiv 1 \mod 4$; quadratic reciprocity yields
$$
\leg 3p = \leg p3,\qquad \leg 5p = \leg p5 \ ;
$$
and as $p \equiv 17 \mod 60$, we have $p \equiv 2 \mod 3, 5$ so the supplement
to quadratic reciprocity gives $-1$ for each.
\msk
Note that since $p$ is not $3$ nor $5$, $u+u^{-1} = 1/2$ cannot be
$\pm 2$.  As $WV$, $WV^{-1}$ and the diagonal matrix $\Delta$
with entries $u,\ u^{-1}$ have equal traces, they are conjugate.
\msk
{\bf Lemma 7.} \ $c(W\gen{V}) = (p+1)/2\ cl(V) + (p-1)\ cl(WV)$.
\msk
{\bf Proof:} \ First note that $|N(W)| = p+1$ and $|N(WV)| = p-1$ by
Lemma 2.  Clearly then $|N(\{WV, WV^{-1}\})|$ is a divisor of
$2(p-1)$.  Hence $N(W\gen{V})$ has order dividing both $p+1$ and
$2(p-1)$, that is, dividing $2$, and hence is $\{\pm I\}$.  It follows
that $cl(W) = cl(V)$ appears with coefficient $|N(W)|/|N(W\gen{V})| =
(p+1)/2$.  Subtracting this from the total number of elements in
$c(W\gen{V})$ leaves
$$
(p-1)p(p+1) = 3\cdot (p-1)p(p+1)/2 - (p+1)/2\cdot p(p-1) \ .
$$
But $cl(WV)$ has $p(p+1)$ elements by Lemma 2, so must appear $p-1$
times.  It follows then that $N(\{WV,WV^{-1}\})$ must have exactly
$2(p-1)$ elements.
\msk
Finally, equation $(*)$ now follows immediately from Lemmas 4--7.


\vfill
\newpage
\bsk
\centerline{\bf References}
\bsk


\alwaysindent=15pt	
\ssk

\ssk
\def\berger{\hbox to 26pt{[Be]\hss}}
\item{\berger} Ruth I\. Berger, 
Ideal class groups, some computations.
J\. Number Theory 50 (1995), no\. 2, 251--260. 

\ssk
\def\boge{\hbox to 26pt{[B\"o]\hss}}
\item{\boge} Sigrid B\"oge, {\it Die $\epsilon$-Invariante von $\text{SL}(2,p)$,}
Arch\. Math\. (Basel) {\bf 46} (1986), no.~4,
299--303. 
\ssk
\def\brown{\hbox to 26pt{[Br]\hss}}
\item{\brown} K\. Brown, Groups with a fixed point free
representation, unpublished lectures, 25 pp\. Spring, 1975.

\ssk
\def\rbrown{\hbox to 26pt{[Br1]\hss}}
\item{\rbrown} Ron Brown, 
Frobenius groups and classical maximal orders, 
Mem\. Amer\. Math\. Soc\., {\bf 151} (2001), no\. 717, viii+110.

\ssk
\def\rbrownh{\hbox to 26pt{[BrH]\hss}}
\item{\rbrownh} Ron Brown and David K\. Harrison, 
Abelian Frobenius kernels and modules over number rings, 
J\. Pure Appl\. Algebra, {\bf 126} (1998), no\. 1-3, 51--86.

\ssk
\def\cgw{\hbox to 26pt{[CGW]\hss}}
\item{\cgw} J\.-L\. Colliot-Th\'el\`ene, R\. M\. Guralnick and R\. Wiegand, 
Multiplicative groups of fields modulo products of subfields, 
J\. Pure Appl\. Algebra, {\bf 106} (1996), no\. 3, 233--262.

\ssk
\def\cgwc{\hbox to 26pt{[CGW2]\hss}}
\item{\cgwc} J\.-L\. Colliot-Th\'el\`ene, R\. M\. Guralnick and R\. Wiegand, 
Erratum: ``Multiplicative groups of fields modulo products of subfields'', 
J\. Pure Appl\. Algebra, {\bf 109} (1996), no\. 3, 331.

\ssk
\def\desmit{\hbox to 26pt{[dS]\hss}}
\item{\desmit} Bart de Smit, 
Primitive elements in integral bases, 
Acta Arith\., {\bf 71} (1995), no\. 2, 159--170.

\ssk
\def\dorn{\hbox to 26pt{[D]\hss}}
\item{\dorn} Larry Dornhoff, {\it Group representation theory: Part
A. Ordinary representation theory. Marcel Dekker, Inc\. New York,
1971.} 
\ssk
\def\gapx{\hbox to 26pt{[GAP]\hss}}
\item{\gapx} GAP - Groups, Algorithms and Programming\hfill\break
https://www.gap-system.org/

\ssk
\def\guwi{\hbox to 26pt{[GW]\hss}}
\item{\guwi} Robert Guralnick and Roger Wiegand, {\it Galois groups and the
             multiplicative structure of field extensions,}
             Trans\. Amer\. Math\. Soc\. {\bf 331} (1992), no.~2, 563--584.

\ssk
\def\hsia{\hbox to 26pt{[HP1]\hss}}
\item{\hsia} J.~S.~Hsia and Roger D.~Peterson, {\it An invariant ideal of a 
             group ring of a finite group, and
             applications,} J.~Algebra {\bf 32} (1974), no.~3,
             576--599. 

\ssk
\def\hsiap{\hbox to 26pt{[HP2]\hss}}
\item{\hsiap}  J.~S.~Hsia and Roger D.~Peterson, 
               {\it An invariant ideal of a group ring of a finite
               group. II,}  Proc.~Amer.~Math.~Soc\. {\bf 51} (1975), no.~2,
               275--281. 

\ssk
\def\kata{\hbox to 26pt{[K]\hss}}
\item{\kata} Toshitaka Kataoka, {\it Some type of ideals of group rings and its
             applications to algebraic number fields,} J\. Faculty of
             Science, Univ\. Tokyo, Sec\. IA, no.~3 {\bf 26} (1979),
             443-452. 

\ssk
\def\rehmone{\hbox to 26pt{[RH]\hss}}
\item{\rehmone} Hans Peter Rehm and Wolfgang Happle,
Zur gruppentheoretischen Absch\"atzung von Idealklassenexponenten
galoisscher Zahlk\"orper durch Exponenten geeigneter Teilk\"orper,
J\. Reine Angew\. Math\., {\bf 268/269} (1974), 434--440.

\ssk
\def\rehmtwo{\hbox to 26pt{[R]\hss}}
\item{\rehmtwo} Hans Peter Rehm,
\"Uber die gruppentheoretische Struktur der Relationen zwischen Relativnormabbildungen in endlichen Galoisschen K\"orpererweiterungen,
J\. Number Theory, {\bf 7} (1975), 49--70.

\ssk
\def\scharlaua{\hbox to 26pt{[Sch 1]\hss}}
\item{\scharlaua} Winfried Scharlau, Induction theorems and the structure
             of the Witt group., Invent\. Math\. {\bf 11} (1970), 37--44.

\ssk
\def\scharlaub{\hbox to 26pt{[Sch 2]\hss}}
\item{\scharlaub} M.~Knebusch and W.~Scharlau, \"Uber das Verhalten der
                 Witt-Gruppe bei galoischen K\"orpererweiterungen,
                 Math\. Ann\. {\bf 193} (1971), 189--196.

\ssk
\def\scharlau{\hbox to 26pt{[Sch]\hss}}
\item{\scharlau} Winfried Scharlau, {\it Eine Invariante endlicher Gruppen,}
        	 Math\. Z\. {\bf 130} (1973),
        	 291--296.
	
\ssk
\def\schur{\hbox to 26pt{[S]\hss}}
\item{\schur} Issai Schur, {\it Untersuchungen \"uber die
                 Darstellungen der endlichen Gruppen durch gebrochene
                 lineare Substitutionen,}  J. Reine Angew\. Math\. {\bf
                 132} (1907), 85--137.  (in [SW], pp\. 198--250)
  		 
\ssk
\def\schurw{\hbox to 26pt{[SW]\hss}}
\item{\schurw} Isaai Schur, {\it Gesammelte Abhandlungen. Band I},
	       eds\. Alfred Brauer and Hans Rohrbach.  Springer,
	       Berlin, 1973. 	       

\ssk
\def\suzuki{\hbox to 26pt{[Su]\hss}}
\item{\suzuki}  Michio Suzuki, {\it On finite groups with cyclic sylow subgroups for all
                odd primes,} American J.~Math., {\bf 77}, no.~4 (1955),
                657-691. 

\ssk
\def\wolf{\hbox to 26pt{[W]\hss}}
\item{\wolf} Joseph A\. Wolf, Spaces of constant curvature. Third edition.
		Boston, Mass\., 1974. Publish or Perish, Inc\. xv+408
                pp\. 
\ssk
\def\young{\hbox to 26pt{[Y]\hss}}
\item{\young} Paul K\. Young, The Scharlau Invariant, Senior Thesis,
                Cornell University, 2003, 22 pp\.
\ssk
\def\zass{\hbox to 26pt{[Z]\hss}}
\item{\zass}
Hans Zassenhaus, \"Uber endliche Fastk\"orper,  Abh\. Math\. Sem\. Univ\. Hamburg,
{\bf 11} (1936), 187--220.  

\bsk

R.~Keith Dennis, Department of Mathematics, Cornell University,
Ithaca, NY 14853\

Email address:\ \ dennis\@math.cornell.edu

\ssk
Paul K.~Young, 964 Hancock Ave \#303
West Hollywood, CA 90069\

  Email address:\ \ paulkyoung01\@gmail.com

\bye

%% file: fixitem.tex
%
\newdimen\alwaysindent
\alwaysindent=20pt
\newbox\anindent
\setbox\anindent=\hbox to\alwaysindent{}
\def\hang{\hangindent\alwaysindent}
\def\textindent#1{\hskip\alwaysindent\llap{#1\enspace\thinspace}\ignorespaces}
\def\item{\par\noindent\hang\textindent}

\def\narrower{\advance\leftskip\alwaysindent
  \advance\rightskip\alwaysindent}
%
%

%